\documentclass[12pt,oneside]{amsart}
\usepackage[english]{babel}
\usepackage{amssymb}
\usepackage{amsthm}
\usepackage{mathtools}

\usepackage[foot]{amsaddr}

\usepackage{amsmath}
\usepackage{a4wide}
\usepackage{nccmath}
\usepackage{esvect}
\usepackage{hyperref}
\hypersetup{
    colorlinks,
    linkcolor=red,
    citecolor=black,
    urlcolor=blue
}
\usepackage{verbatim}
\usepackage[dvips]{graphicx}
\usepackage{t1enc}
\usepackage{color}
\usepackage{geometry}
\usepackage{graphicx}
\usepackage{color,soul}
\usepackage[parfill]{parskip}   
\usepackage[numbers,sort&compress]{natbib}
\usepackage{tikz}
\usepackage{amsmath}
\usepackage{caption}
\usepackage{subcaption}
\usepackage{mathtools}
\usepackage{subfiles}

\usetikzlibrary{graphs}
\usetikzlibrary{graphs.standard}

\makeatletter
\def\subsubsection{\@startsection{subsubsection}{3}%
  \z@{.5\linespacing\@plus.7\linespacing}{.1\linespacing}%
  {\normalfont\itshape}}
\makeatother

\newtheorem{thm}{Theorem}[section]
\newtheorem{lem}[thm]{Lemma}
\newtheorem{lemma}[thm]{Lemma}

\newtheorem{prop}[thm]{Proposition}

\newcommand{\EE}{\mathbb{E}}

\newcommand{\PP}{\mathbb P}

\newcommand{\ignore}[1]{}

\newcommand{\rh}[1]{{\color{orange}{#1}}}
\newcommand{\pvh}[1]{{\color{red}{#1}}}

\begin{document}

\author{Rebekah Herrman}
\address[Rebekah Herrman]{Department of Industrial and Systems Engineering, The University of Tennessee, Knoxville, TN.}
\email[Rebekah Herrman]{rherrma2@utk.edu}

\author{Peter van Hintum}
\address[Peter van Hintum]{Department of Pure Maths and Mathematical Statistics, University of Cambridge, UK}
\email[Peter van Hintum]{pllv2@cam.ac.uk}

\author{Stephen G. Z. Smith}
\email[Stephen G. Z. Smith]{sgsmith1@memphis.edu}

\title[Capture times in the  Bridge-burning Cops and Robbers game] 
{Capture times in the  Bridge-burning Cops and Robbers game}

\linespread{1.3}
\pagestyle{plain}

\begin{abstract}
In this paper, we consider a variant of the cops and robbers game on a graph, introduced by Kinnersley and Peterson, in which every time the robber uses an edge, it is removed from the graph, known as bridge-burning cops and robbers. In particular, we study the maximum time it takes the cops to capture the robber.
\end{abstract}
\maketitle 

\section{Introduction}


\textit{Cops and Robbers} is a well-studied game on a graph $G$ with two players, the \emph{cops} and the \emph{robber}. The game begins with the cops choosing their starting vertices, followed by the robber selecting his. The cops and robber then alternate turns moving from their current vertex to an adjacent one or choosing not to move. A \emph{round} consists of a cop turn and the subsequent robber turn. The cops win, or capture the robber, if a cop and the robber occupy the same vertex, whereas the robber wins if he manages to avoid capture.

Several variants of the game have been introduced over the years, including those where the robbers can move more quickly than the cops \cite{balister, frieze}, where the cops have imperfect information \cite{clarke}, and where only one cop may move during the cops' turn \cite{das}. Recently, Kinnersley and Peterson \cite{kinnersley} introduced the variant \textit{bridge-burning cops and robbers}. In the bridge-burning version, each time the robber moves from a vertex $u$ to a vertex $v$, the edge $uv$ is erased from the graph. Using the notation introduced in \cite{kinnersley},  $c_b(G)$ is defined to be the \emph{bridge burning cop number}, which is the minimum number of cops required to catch the robber on the graph $G$ in the bridge-burning game. Kinnersley and Peterson studied the game on numerous graphs including paths $P_n$, cycles $C_n$, complete bipartite graphs $K_{m,n}$, hypercubes $Q_n$, and two dimensional finite grids $G_{m,n}$ \cite{kinnersley}. 

A related notion to $c_b(G)$ is the \emph{capture time} of $G$, denoted $capt_b(G)$. The bridge-burning capture time is the minimum number of rounds it takes for the cop to capture the robber in the bridge-burning variant. The capture time of the original cops and robbers game was introduced in 2009 by Bonato, Golovach, Hahn, and Kratochv\'il \cite{bonato} and has been studied on trees \cite{yang} and planar graphs \cite{pisantechakool} among various other classes of graphs \cite{gavenvciak2010cop, mehrabian2011capture, bonato2013capture, kinnersley2018bounds}. Counting the number of possible configurations shows quickly that for graphs with cop number $c(G)\leq k$, we have capture time $\text{capt}(G)=O(n^{k+1})$, where $n$ is the number of vertices in $G$. In their original paper \cite{bonato}, Bonato, Golovach, Hahn, and Kratochv\'il showed that for $c(G)=1$, this can be improved to $\text{capt}(G)=O(n)$, which is tight as shown by the path graph $P_n$. Perhaps surprisingly, Brandt, Emek, and Uitto \cite{brandt} showed that for $c(G)\geq2$, the trivial upper bound is actually tight, i.e. there exist graphs $G$, with $\text{capt}(G)=\Omega(n^{c(G)+1})$, indicating a qualitative difference between graphs $G$ with $c(G)=1$ and $c(G)\geq 2$.

Returning our attention to the bridge-burning capture time, Kinnersley and Peterson \cite{kinnersley} showed that if $c_b(G)=1$, then $capt_b(G)= O(n^3)$ and conjectured that there exists a graph $G$ such that $c_b(G)=1$ and $capt_b(G)= \Omega(n^3)$. We generalise their result by showing that $capt_b(G)=O(n^{c_b(G)+2})$ and prove the matching lower bound analogous to the one in their conjecture for $c_b(G)\geq 3$.
 \begin{thm}\label{mainthm}
 There exists a universal constant $C>0$ such that the following holds. For every $k\geq 3$ and $n$ sufficiently large, there exists a graph $G_n$ such that $v(G_n)=n$, $c_b(G_n)=k$, and 
 $$C\ \frac{n^{k+2}}{k^{k+2}}\leq capt_b(G_n).$$
 \end{thm}
 
In fact, in Proposition \ref{upperbound} we show that for all $G$ on $n$ vertices $capt_b(G)\leq  \frac{(2n)^{c_b(G)+2}}{c_b(G)!}$, which shows that  the asymptotics in terms of $c_b(G)$ are almost tight.

In Section \ref{capturetimes}, we present some preliminary results on catching times in the bridge-burning game and in Section \ref{proofsection} we prove Theorem \ref{mainthm}.

As usual, we write $[k]=\{1,\dots,k\}$. The asymptotics in the paper are with respect to the number of vertices $n$, assuming fixed burning-bridge cop number, unless explicitly stated to be otherwise.
 


  \section{Capture times} \label{capturetimes}

In this paper, we will show that the graph $G$ on $n$ vertices with cop number $k\geq 3$ which maximizes the capture time satisfies
$$C\  \frac{n^{k+2}}{k^{k+2}}\leq \max \{capt_b(G):c_b(G)=k, v(G)=n\}\leq C'\ \frac{(2n)^{k+2}}{k!}$$
for some universal constants $C,C'>0$.
 
 First, we show that the capture time of $K_{n,n}$ is $\Theta(n^2)$. Our proof significantly simplifies the proof given in \cite[Theorem 5.2]{kinnersley} to demonstrate that there are graphs with $c(G)=1$ and $capt_b(G)=\Omega(n^2)$. In order to prove the result, we need the following slight strengthening of a theorem from Kinnersley and Peterson \cite[Theorem 2.2]{kinnersley}. 
 \begin{lem}\label{clique}
 If $\exists \; X\subset V(G)$, such that $G[X]$ is a clique and $X\cup\Gamma(X)=V(G)$, then $c_b(G)=1$ and $capt_b(G)=O(n^2)$, where $\Gamma(X)$ is the neighbourhood of $X$.
\end{lem}
\begin{proof}
 Place the cop on any vertex in $X$. Subsequently, always move the cop to a vertex in $X$ adjacent to the position of the robber.
Note that the robber can never move onto a vertex in $X$ and, thus, can never remove an edge incident to $X$. Hence, $X\cup\Gamma(X)=V(G)$ remains constant throughout the game. After each round, the cop is adjacent to the robber, so the robber must move in every round. Given that the robber removes one edge in every round, eventually he must move into $X$, as all the other possible edges have been removed. As there are $O(n^2)$ edges, this must happen within $O(n^2)$ moves.
\end{proof}

This lemma provides the cop number and an upper bound in the following proposition.
 \begin{prop} $K_{n,n}$ has capture time $\Theta(n^2)$\end{prop}
 
 \begin{proof}
 As any two adjacent vertices in $K_{n,n}$ satisfy the conditions in Lemma \ref{clique},  we find that $c_b(K_{n,n})=1$ and $capt_{b}(G)=O(n^2)$.
 
 On the other hand, we consider the following strategy for the robber to delay capture. First, we find an Euler cycle of $K_{\lfloor\frac{n}{2}\rfloor,\lfloor\frac{n}{2}\rfloor}$ (or $K_{\lfloor \frac{n}{2}\rfloor-1,\lfloor \frac{n}{2}\rfloor-1}$ if $\lfloor \frac{n}{2}\rfloor$ is odd). Next, we traverse the following route through $K_{n,n}$; to each vertex in $K_{\lfloor\frac{n}{2}\rfloor,\lfloor\frac{n}{2}\rfloor}$, assign a distinct pair of vertices in $K_{n,n}$ such that the pairs of vertices from the same part of $K_{\lfloor\frac{n}{2}\rfloor,\lfloor\frac{n}{2}\rfloor}$ are in the same part of $K_{n,n}$. Now, the robber follows the Euler cycle through $K_{n,n}$ in the sense that every time he is forced to move, he goes to an element in the corresponding pair which is available. As the cop is only able to occupy one vertex of a given pair, there is no way for the cop to obstruct the robber's path. This route has length $\Omega(n^2)$.
 \end{proof}
 
 In fact Lemma \ref{clique} implies the following result for random graphs $\mathcal{G}(n,p)$.
 
\begin{prop}
 Let $G\in \mathcal{G}(n,p)$. Then w.h.p. $c_b(G)=1$ and $capt_b(G)=O(n^2)$.
\end{prop}
\begin{proof}
 By Lemma \ref{clique}, it suffices to show that $G$ contains a dominating clique with high probability. This follows from a second moment argument included in Lemma \ref{randomlem} in the Appendix.
\end{proof}

For general graphs, we find the following generalization of a result from \cite[Theorem 5.1]{kinnersley} which showed this proposition in the case $c_b(G)=1$. 
 \begin{prop}\label{upperbound}
Let $G$ be a graph on $n$ vertices, then $capt_b(G)\leq\frac{(2n)^{c_b(G)+2}}{c_b(G)!}$. \end{prop}
 
 \begin{proof}
 Note that as the robber removes an edge with every move, the robber can make at most $e(G)\leq \binom{n}{2}$ moves before getting caught. Between two moves of the robber, the cops move around. Without the robber moving, there is no point in the cops returning twice to the exact same configuration. As there are at most $\binom{n+c_b(G)-1}{c_b(G)}\leq\frac{(2n)^{c_b(G)}}{c_b(G)!}$ configurations of the cops on the vertices, it can take at most $\frac{(2n)^{c_b(G)+2}}{c_b(G)!}$ moves before the robber is caught.
 \end{proof}

The remainder of the paper is dedicated to proving Theorem \ref{mainthm}.

 \section{Proof of Theorem \ref{mainthm}}\label{proofsection}

We first prove the result for $k=3$ and then extend the construction to larger $k$. We claim the following graph $G$ has $c_b(G)=3$ and $capt_b(G)=\Theta(n^5)$.
Let the vertex set of $G$ be the following union of sets;

\begin{align*}
   V(G)=\{p_i,&q_i:i\in[3n]\}\cup\{x_1,x_2\}\cup \{d_{x},h_x\}\cup X \cup Y\\
   &\cup \{d_{X,1},d_{X,2},d_{Y,1},d_{Y,2},h_{X,1},h_{X,2},h_{Y,1},h_{Y,2}\}\\
   &\cup \{a_i:i\in[3n]\}\cup\{d_a,h_a\}\cup \{d_{a,v,1},d_{a,2},h_{a,1},h_{a,2}\}\\
   &\cup \{b_i:i\in[3n]\}\cup \{d_{b,1},d_{b,2},h_{b,1},h_{b,2}\}\\
   &\cup \{d_{a,i,1},d_{a,i,2},h_{a,i,1},h_{a,i,2}: i\in[3]\}\cup \{d_{b,i,1},d_{b,i,2},h_{b,i,1},h_{b,i,2}: i\in[3]\}
\end{align*}
 
 with $|X|=|Y|=3n$. Let the edge set of $G$ be the following union of sets;
\begin{align*}
    E(G)&=\{p_ip_{i+1},q_iq_{i+1}:i\in[3n-1]\}\cup\{p_1x_1,q_1x_1,p_nx_2,q_nx_2\}\\
    &\cup\{p_id_x,q_id_x,x_1d_x,x_2d_x:i\in[3n]\}\\
    &\cup \{x_1v:v\in X\}\cup \{x_2v: v\in Y\}\cup \{uv:u\in X, v\in Y\}\\
    &\cup \{p_id,q_id,x_1d,x_2d:i\in[3n],d\in\{d_{X,1},d_{X,2},d_{Y,1},d_{Y,2}\}\}\\
    &\cup\{vd_{X,1},vd_{X,2}:v\in X\}\cup \{vd_{Y,1},vd_{Y,3}:v\in Y\}\\
    &\cup \{a_ia_{i+1}: i\in[3n]\}\cup \{a_id_a:i\in[3n]\}\cup \{a_id_{a,1},a_id_{a,2}: i\in[3n]\}\\
    &\cup \{vd_{a,1},vd_{a,2}:v\in X\cup Y\}\cup \{a_ix_1: i\in[3n]\}\\
    &\cup \{b_ib_{i+1}: i\in[3n]\}\cup \{b_id_b:i\in[3n]\}\cup \{b_id_{b,1},b_id_{b,2}: i\in[3n]\}\\
    &\cup \{vd_{b,1},vd_{b,2}:v\in X\cup Y\}\cup \{b_ix_2:i\in[3n]\}\\
    &\cup \{vd_{a,i,1},vd_{a,i,2}:v\in X\cup Y, i\in[3]\}\\
    &\cup \{x_1d_{a,i,1},x_1d_{a,i,2},x_2d_{a,i,1},x_2d_{a,i,2}:i\in[3]\}\\
    &\cup \{p_jd_{a,i,1},p_jd_{a,i,2},q_jd_{a,i,1},q_jd_{a,i,2}: j\not\equiv i\mod 3\}\\
    &\cup \{a_{(i-1)n+j}d_{a,i,1},a_{(i-1)n+j}d_{a,i,2}:  j\in[n], i\in[3]\}\\
    &\cup \{vd_{b,i,1},vd_{b,i,2}:v\in X\cup Y, i\in[3]\}\cup \{a_jd_{b,i,1},a_jd_{b,i,2}: j\not\equiv i\mod 3\}\\
    &\cup \{b_{(i-1)n+j}d_{b,i,1},b_{(i-1)n+j}d_{b,i,2}:  j\in[n], i\in[3]\}\\
    &\cup \{h_id_i:\text{all $i$, such that $d_i\in V(G)$}\}
\end{align*}

\begin{figure}
\centering
\includegraphics[scale=.56]{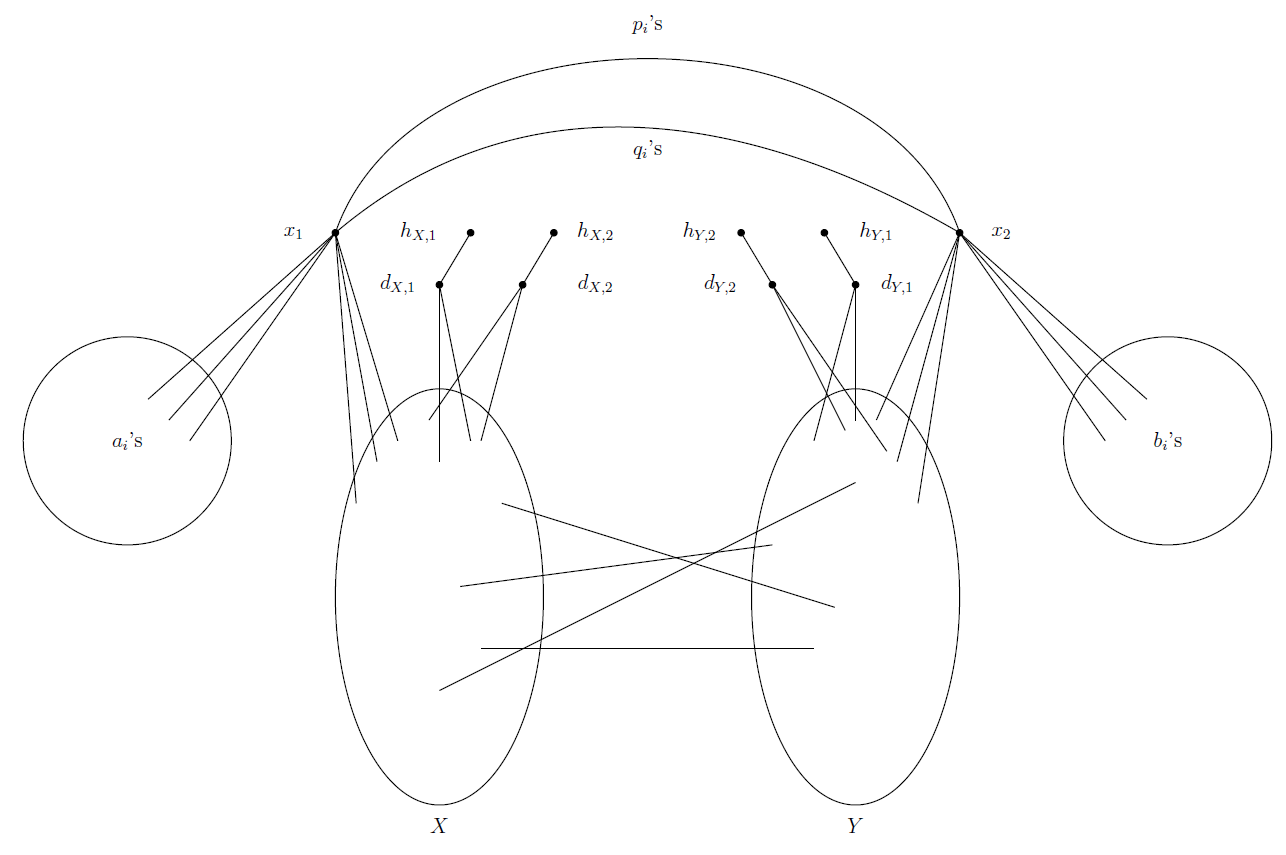}
\caption{The graph $G$ described in the proof of Theorem \ref{mainthm}. Most of the doors and holes are omitted though all the other vertices and edges are displayed.}
\label{main}
\end{figure}

For an illustration of $G$, see Figure \ref{main}.

 In the remainder of the paper we shall use the notation $\{a_i\}_i$ to denote the set of $a_i$'s with $i$ ranging over all possible values in the given context, so e.g. $\{p_i,q_i,x_i\}_i=\{p_i,q_i,x_j:i\in[3n],j=1,2\}$.

This graph $G$ consists of three cycles, $\{a_i\}_i$, $\{b_i\}_i$ and $\{x_i,p_i,q_i\}_i$, a complete bipartite graph on the sets $X$ and $Y$ and a great number of \emph{doors}, $d_i$'s, and \emph{holes}, $h_i$'s. Holes are vertices with degree one and doors are their unique neighbours. 

One of the cycles, viz $\{x_i,p_i,q_i\}_i$, contains two special vertices, $x_1$ and $x_2$, each of which is complete to one of the parts of the bipartite graph and to one of the cycles. In particular, we have $\{a_i\}_i\cup X\subset \Gamma(x_1)$ and $\{b_i\}_i\cup Y\subset \Gamma(x_2)$. In fact, these are the only edges between vertices that are not doors or holes.

The doors and holes restrict the freedom of the cops; if the robber manages to move to a door that is not adjacent to a cop, he will move to the corresponding hole in the next move, disconnecting himself from the rest of the graph and thus winning the game.

Clearly, $G$ has $O(n)$ vertices. We first aim to establish that $c_b(G)=3$, starting with the lower bound.

\begin{lemma}\label{lowerbound3}\label{standardposition}
$c_b(G)\geq 3$. Moreover, if $c_b(G)=3$, then one cop starts in $\{a_i\}_i$, one cop starts in $\{b_i\}_i$  and one cop starts in $\{x_i,p_i,q_i\}_i$.
\end{lemma}
\begin{proof}
To see that $c_b(G)\geq 3$, note that we initially need a cop next to, or on, every door. In particular, doors $d_a,d_b$ and $d_{x}$. Since $\Gamma(d_a)=\{a_i\}_i\cup\{h_a\}$,  $\Gamma(d_b)=\{b_i\}_i\cup\{h_b\}$ and $\Gamma(d_x)=\{x_i,p_i,q_i\}_i\cup\{h_x\}$, there is no vertex next to or on more than one of these, so we need at least three cops.
If $c_b(G)=3$, then evidently one cop must start in each of $\{d_a\}\cup\Gamma(d_a)$, $\{d_b\}\cup\Gamma(d_b)$, and $\{d_x\}\cup\Gamma(d_x)$. If one of these cops doesn't start in the corresponding cycle ($\{a_i\}_i$,$\{b_i\}_i$ and $\{x_i,p_i,q_i\}_i$ respectively), then no cop is adjacent to some other door ($d_{a,1}$, $d_{b,1}$ or $d_{X_1}$ respectively).
\end{proof}

To see that three cops suffice to catch the robber, consider the following strategy for the cops. Start one cop on $a_1$, say Alex, one on $b_1$, say Blake, and one on $x_1$, say Charlie. We will refer to this starting position as the \emph{standard position}. Note that these three vertices cover all the doors. Each of the cops will stay on their respective cycles unless the robber moves onto a vertex adjacent to them, in which case they catch him. 



\begin{figure}
 \centering
 \includegraphics[scale=1]{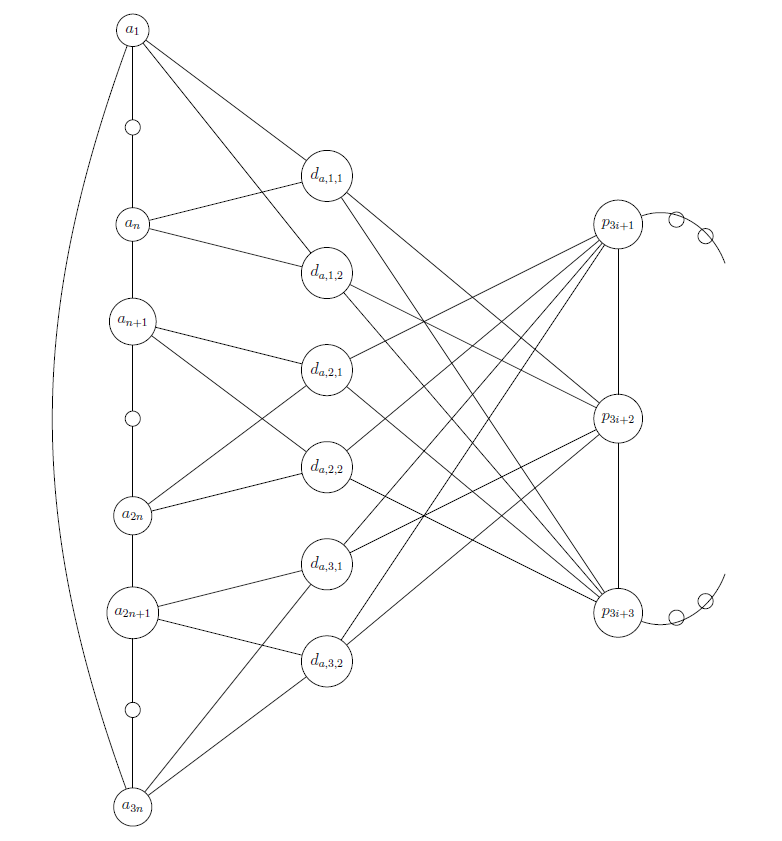}
\caption{The graph described in the proof of Lemma \ref{flexguarding}. The central vertices are doors, and the other vertices form cycles $\{a_i\}_i$ and $\{p_i,q_i,x_i\}$ patrolled by Alex and Charlie respectively which watch the doors.}
\label{B}
\end{figure}
\begin{lemma}\label{flexguarding}\label{manymoves}
If the cops start in the standard position, then every cop can reach any vertex in their cycle, while guarding all doors at each of the intermediate steps. \\
Moreover, if the robber starts and remains in $X\cup Y$ and the cops start in standard position and remain in their cycles always guarding all the doors, then it takes Charlie $\Omega(n^3)$ moves to get from $x_1$ to $x_2$ and from $x_2$ to $x_1$.
\end{lemma}
\begin{proof}
 We will show that every cop can move to a neighbouring vertex in at most $O(n^2)$ steps. Recall that each cycle has diameter $O(n)$.
 
We first consider Blake's moves. Blake can move freely between the vertices in $\{b_{kn+j}:j\in[n]\}$ for any fixed $k\in\{0,1,2\}$, as each of these vertices has the same neighbourhood outside $\{b_i\}_i$. When changing $k$, Blake's neighbourhood in $\{d_{b,i,j}\}_{i,j}$ changes, so in order to keep guarding all the doors, Alex must move in parallel to cover Blake's old neighbourhood (cf. Figure \ref{B}). This, in turn, affects Alex's neighbourhood in $\{d_{a,i,j}\}_{i,j}$, which would then have to be compensated by Charlie. Thus, Blake can move anywhere in $\{b_i\}_i$ in $O(n)$ steps.

For Alex, the concerns are very similar. Two adjacent vertices in $\{a_i\}_i$ have different neighbourhoods in $\{d_{b,i,j}\}_{i,j}$, so for every consecutive step Alex takes, Blake has to take $n$ steps. Hence, Alex can move anywhere in $O(n^2)$ steps. Moreover, to move to a vertex at distance $\Omega(n)$ in the cycle $\{a_i\}$ takes $\Omega(n^2)$ moves.

Finally, every move by Charlie requires $O(n)$ steps of Alex, which in turn requires $O(n^2)$ steps by Blake. Hence, Charlie can move anywhere in $O(n^3)$. Moreover, to move from $x_1$ to $x_2$ and back takes $\Omega(n^3)$ moves.
\end{proof}

We need to exclude the case that the robber doesn't start in $X\cup Y$.

\begin{lemma}\label{outsideXY}
If the cops start in the standard position and the robber starts on any vertex that is not in $X\cup Y$, then the robber is caught in $O(n^3)$ moves.
\end{lemma}
\begin{proof}
If the robber starts on a door or hole or $x_1$ or $x_2$, then the cops can immediately catch or corner him.

Alternatively, suppose the robber starts in one of the cycles. If the cops stay in their cycles, they can move along the cycles while guarding all the doors as shown in the previous lemmas. This implies that the robber cannot leave the cycle he starts in. It is easy to catch a robber on a cycle in $O(n)$ moves. Every step by the cop can require up to $O(n^2)$ moves by the other cops, so the cops need $O(n^3)$ moves to catch the robber.
\end{proof}

Now that we may assume the robber starts in $X\cup Y$, we are ready to show that Alex, Blake and Charlie will succeed in catching the robber.

\begin{lemma}
$c_b(G)=3$
\end{lemma}
\begin{proof}
Lemma \ref{lowerbound3} shows we need at least three cops, so we only need provide a bound from above.

Consider the cops starting in the standard position. If the robber starts outside $X\cup Y$, then the cops can catch the robber according to Lemma \ref{outsideXY}. Hence, we may assume the robber starts in $X\cup Y$.

The cops will move in such a way that all doors are guarded at all times. Moreover, Alex and Blake will stay on $\{a_i\}_i$ and  $\{b_i\}_i$ respectively at all times. Hence, if at any point the robber leaves the set $X\cup Y$, either to a door or to one of $x_1,x_2$, then the cops can immediately seize him. Hence, the robber has to stay inside $X\cup Y$.

Finally, to show that the cops can actually capture the robber, it suffices to show that they can force the robber to keep moving, as he can make at most $|X|\cdot |Y|$ moves staying on the vertices of $X\cup Y$. To this end, Charlie will move between $x_1$ and $x_2$, which by Lemma \ref{flexguarding} is possible while ensuring the cops guard all the doors at every intermediate step. As $x_1$ is complete to $X$ and $x_2$ is complete to $Y$, this forces the robber to keep moving. Hence, the cops will eventually capture the robber.
\end{proof}


To find the lower bound on the capture time, we need to be sure that the cops do not have a better strategy.

\begin{figure}
 \centering
 \includegraphics[scale=.6]{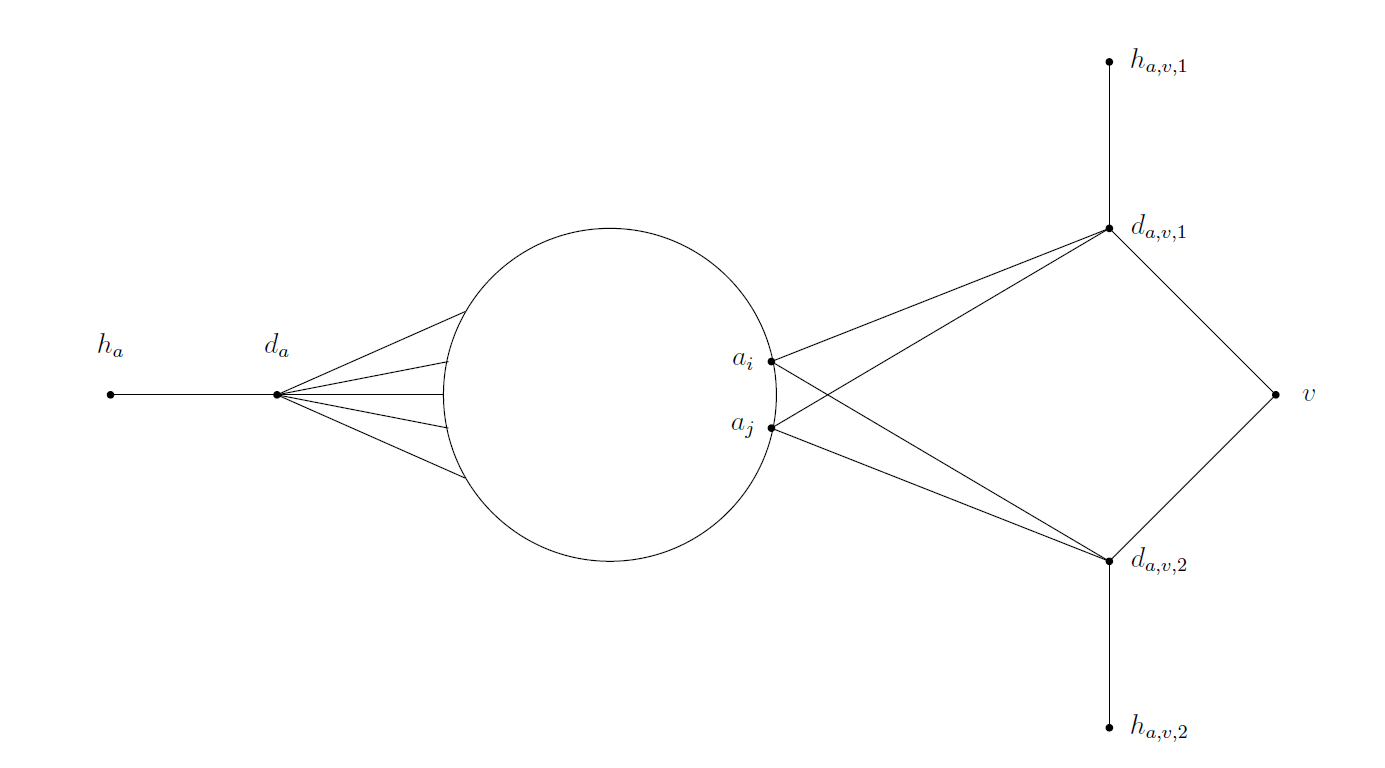}
\caption{The graph described in Lemma \ref{coprestriction}, where $v$ is a vertex in $X\cup Y$.}
\label{A}
\end{figure}

\begin{lemma} \label{coprestriction}
If the cops start in the standard position and the robber starts in $X\cup Y$, then the cops need to stay on their respective cycles for as long as the robber stays in $X\cup Y$, unless they can capture the robber directly.
\end{lemma}
\begin{proof}
For Charlie, let the robber be on $v\in X$, without loss of generality. If Charlie leaves the cycle still guarding $d_{X,1}$ and $d_{X,2}$, then Charlie must have moved into $X$. However, that would imply the cop was previously on $x_1$, so Charlie could have caught the robber immediately. Hence, Charlie cannot leave the cycle without allowing the robber to escape.

For Alex (resp. Blake), note that if the robber is on $v\in X\cup Y$, then leaving their cycles would mean leaving either $d_{a,v,1}$ or $d_{a,v,2}$ (resp. $d_{b,v,1}$ or $d_{b,v,2}$) unguarded, providing an escape route for the robber, as seen in Figure \ref{A}.
\end{proof}
\pagebreak

\begin{lemma}
$\text{capt}_b(G)=\Omega(n^5)$
\end{lemma}
\begin{proof}
By Lemma \ref{standardposition}, the cops must start in standard position or equivalent. 

The robber will follow the following strategy. He will fix a walk of length $\Omega(n^2)$ through the induced complete bipartite graph on vertex set $X\cup Y$, which he can trivially do. He will proceed to follow this walk as slowly as possible, i.e. only proceeding to the next vertex when a cop is adjacent to him.

The robber will only move through $X\cup Y$, so by Lemma \ref{coprestriction} the cops are confined to their cycles. Only Charlie can be adjacent to $X\cup Y$ without leaving his cycle, so it is up to Charlie to walk up and down between $x_1$ and $x_2$ to force the robber to move. By Lemma \ref{manymoves}, it thus takes the cops $\Omega(n^3)$ moves to make the robber move once. Hence, the robber manages to stay out of the cops hands for $\Omega(n^5)$ moves.
\end{proof}

The construction slowing down Charlie can be extended in a natural way to higher cop numbers. Consider the following construction for cop number $k$. For $k\geq 3$, let $G_k$ be the graph constructed as follows.
\begin{align*}
V(G_k)=V(G)&\cup \{u^j_i:i\in[3n], j\in[k-3]\}\cup \{d_{u^j},h_{u^j}:j\in[k-3]\}\\
&\cup \{d_{u^j,i,l},h_{u^j,i,l}:i\in[3],l\in[2], j\in[k-3]\}\\
E(G_k)=E(G)&\cup\{u^j_{i}u^j_{i+1}:i\in[3n],j\in[k-3]\}\\
&\cup\{u^j_id_{u^j}:i\in[3n],j\in[k-3]\}\\
&\cup \{vd_{u^j,i,l}:v\in X\cup Y, i\in[3], l\in[2],j\in[k-3]\}\\
&\cup \{u^{j-1}_ld_{u^j,i,k}: l\not\equiv i\mod 3,k\in[2],j\in[k-3]\}\\
&\cup \{u^j_{(i-1)n+l}d_{u^j,i,k}:  l\in[n], i\in[3],k\in[2],j\in[k-3]\}
\end{align*}
where $u^{0}_i=b_i$.
The $\{u_i^j\}_i$ form cycles, which are similar to cycles $\{a_i\}_i$ and $\{b_i\}_i$. The doors $\{d_{u^j,i,l}\}_{i,l}$ are connected to respective cycles in the same fashion $\{a_i\}_i$ and $\{b_i\}_i$ are connected to the doors $\{d_{b,i,l}\}_{i,l}$.

\begin{prop}
$c_b(G_k)=k$ and $capt_b(G_k)\geq C v(G_k)^{k+2}k^{-(k+2)}$ for some universal constant $C$.
\end{prop}
\begin{proof}[Sketch of proof]
As in Lemma \ref{standardposition}, each of the doors $d_x,d_a,d_b$ and $d_{u^j}$ (with $j\in[k-3]$) must be guarded initially, so $c_b(G_k)\geq k$. Moreover, if $c_b(G_k)=k$, then one cop must start in each of the cycles; $\{p_i,q_i,x_i\}_i$, $\{a_i\}_i$, $\{b_i\}_i=\{u^0_i\}_i$ and $\{u^j_i\}_i$ for $j\in[k-3]$. As in lemma \ref{outsideXY}, if the robber starts in one of the cycles, he will be captured quickly. If the robber starts in the bipartite graph $X\cup Y$, the cops can prevent him from leaving it. Moreover, by Charlie moving between $x_1$ and $x_2$ the robber can be forced to use up all the edges between $X$ and $Y$ and thus be forced out of the bipartite graph, leading to his immediate capture. Hence, $c_b(G_k) \leq k$.

The robber will follow the same strategy as before; planning out a  Eulerian walk (assume for convenience that $n$ is even) through complete bipartite graph $X\cup Y$ and only proceeding through the walk when Charlie is directly adjacent to him. Note that this walk has length $(3n)^2$. As in Lemma \ref{coprestriction}, the cops are restricted to their cycles as long as the robber stays in $X\cup Y$. As in Lemma \ref{manymoves}, for Charlie to move once from $x_1$ to $x_2$ and back, the cops must make $2(3n)^{k}+o(n)$ moves. Thus, the robber can avoid the cops for at least $(3n)^{k+2}$ rounds.

Note that the graph has $(k+2)3n+18k$ vertices. Hence, 
\begin{align*}
capt_b(G_k)&\geq(3n)^{k+2}\\
&=\left(\frac{v(G_k)-18k}{k+2}\right)^{k+2}\\
&\geq\left(\frac{1-\frac{6}{n}}{1+2/k}\right)^{k+2}\left(\frac{v(G_k)}{k}\right)^{k+2} \\
&\geq C\left(\frac{v(G_k)}{k}\right)^{k+2}
\end{align*}
for some constant $C$.
\end{proof}
This completes the sketch of the proof of Theorem \ref{mainthm}.

\section{Concluding remarks}
The natural question remaining is the asymptotic maximal capture time for $c_b(G)=1,2$. Kinnersley and Peterson \cite{kinnersley} conjectured that there exists an $n$-vertex graph, $G$,  with $c_b(G) = 1$ and $capt_b(G) = \Omega(n^3)$, which we leave open.

Additionally, we are interested in the exact asymptotics in terms of $c_b(G)$. The results in this paper show the function to lie between $\frac{1}{c_b(G)^{c_b(G)+2}}$ and $\frac{2^{c_b(G)}}{c_b(G)!}$. We expect the correct answer to be $\frac{1}{c_b(G)!}$.


 \section*{Acknowledgements}
 We would like to thank our supervisor Professor B\'ela Bollob\'as for his continuing support and comments on an earlier version of this manuscript. The first author would like to thank NSF grants Applications of Probabilistic Combinatorial Methods, DMS-1855745, and Probabilistic and Extremal Combinatorics, DMS-1600742. The second author would like to thank the ESPRC for funding under PhD award 1951104 and the Cambridge trust for the Cambridge European Scholarship 10422838.


\bibliographystyle{abbrv}
\bibliography{references}

\appendix
\section{$G\in\mathcal{G}(n,p)$ has $c_b(G)=1$ whp}

By Lemma \ref{clique}, it suffices to show that whp $G\in\mathcal{G}(n,p)$ contains a dominating clique. We shall abbreviate $\log_{\frac{1}{1-p}}(x)$ to $\log(x)$.

 \begin{lem}\label{randomlem}
  Let $G\in\mathcal{G}(n,p)$ with $p\in(0,1]$ constant, then with high probability, $\exists X\subset V(G)$ such that $G[X]$ is a complete graph and $X\cup \Gamma(X)=V(G)$.
 \end{lem}
 \begin{proof}
 We use a second moment argument to show the result. Fix some small $\epsilon\in(0,\frac12)$ and let $k=(1+\epsilon)\log(n)$.
 
 Let $S$ be the number of sets $X\subset V(G)$ such that $|X|=k$, $X$ induces a clique and $X\cup \Gamma(X)=V(G)$.
  Note that the events that $X$ is a clique and that $X\cup \Gamma(X)=V(G)$ are dependent on disjoint edges.
\begin{align*}
\EE[S]&=\binom{n}{k}p^{\binom{k}{2}}\left(1-\left(1-p\right)^k\right)^{n-k}
\end{align*}
To compute the second moment of $S$, let $A$ and $B$ be two $k$-sets of vertices. We will use the law of total expectation to condition on the size of $A\cap B$. Note that the probability that $A$ and $B$ both satisfy the conditions is at most the probability that they are both independent sets.
\begin{align*}
\EE[S^2]&\leq\sum_{i=0}^k \binom{n}{k-i,i,k-i} p^{2\binom{k}{2}-\binom{i}{2}}\\
&\leq\left(\binom{n}{k}p^{\binom{k}{2}}\right)^2 \left[1+\sum_{i=1}^k \frac{\binom{n}{k-i,i,k-i}}{\binom{n}{k}^2} p^{-\binom{i}{2}}\right]
\end{align*}
Each of these last terms is bounded as:
$$\frac{\binom{n}{k-i,i,k-i}}{\binom{n}{k}^2} p^{-\binom{i}{2}}\leq \frac{k^{2i} p^{-\binom{i}{2}}}{n^i},$$
so for the entire sum we find 
\begin{align*}
\sum_{i=1}^k \frac{\binom{n}{k-i,i,k-i}}{\binom{n}{k}^2} p^{-\binom{i}{2}}\leq \max_{i\in[k]} \left\{ \frac{k^{2i+1} p^{-\binom{i}{2}}}{n^i}\right\}=o(1),
\end{align*}
and thus
$$\EE[S^2]\leq\left(\binom{n}{k}p^{\binom{k}{2}}\right)^2 (1+o(1)).$$

Now we find by Chebyshev's inequality:
\begin{align*}
\PP(S>0)
&\geq \frac{\left[\binom{n}{k} p^{\binom{k}{2}}\left(1-\left(1-p\right)^k\right)^{n-k}\right]^2}{\left(\binom{n}{k}p^{\binom{k}{2}}\right)^2 (1+o(1))}
\to 1
\end{align*}
Hence, the probability that there is a dominating clique tends to one.
\end{proof}

\end{document}